\documentclass[12pt]{amsart}

\newtheorem{theorem}{Theorem}[section]
\newtheorem{lemma}[theorem]{Lemma}

\numberwithin{equation}{section}

\begin{document}

\baselineskip=16pt

\title{Holomorphic Cartan geometries and Calabi--Yau manifolds}

\author[I. Biswas]{Indranil Biswas}

\address{School of Mathematics, Tata Institute of Fundamental
Research, Homi Bhabha Road, Bombay 400005, India}

\email{indranil@math.tifr.res.in}

\author[B. McKay]{Benjamin McKay}

\address{School of Mathematical Sciences,
University College Cork, National University of Ireland}

\email{B.McKay@ucc.ie}

\subjclass[2000]{}

\date{}

\begin{abstract}
We prove that the only Calabi--Yau projective manifolds
which bear holomorphic Cartan geometries are 
precisely the abelian varieties.
\\

\noindent
\textsc{R\'esum\'e.} Nous d\'emontrons que les seules
vari\'et\'es projectives de Calabi--Yau qui poss\`edent
des g\'eom\'etrie holomorphes de Cartan sont les vari\'et\'es
ab\'eliennes.
\end{abstract}

\maketitle

\section{Introduction} 

Let $G$ be a complex Lie group and $H\,\subset\, G$
a closed complex subgroup. The Lie algebras of
$G$ and $H$ will be denoted by $\mathfrak g$ and $\mathfrak h$
respectively.

Let $M$ be a connected complex manifold. Let
$$
E_H\, \longrightarrow\, M
$$
be a holomorphic principal right $H$--bundle. For each
point $h\in\, H$, let
\begin{equation}\label{e2}
\tau_h\, :\, E_H\, \longrightarrow\, E_H
\end{equation}
be the translation by the action of $h$.
For any $v\, \in\, \mathfrak h$, let $\zeta_v$ be the
holomorphic vector field on $E_H$ 
that assigns to $z\, \in\,
E_H$ the tangent vector defined by the analytic curve
$$
t\, \longmapsto\, z\exp(tv)\, .
$$

A \textit{holomorphic Cartan geometry} on $M$ of type $G/H$
is a holomorphic principal $H$--bundle
\begin{equation}\label{e0}
E_H\, \longrightarrow\, M
\end{equation}
together with a holomorphic 1--form on
$E_H$ with values in $\mathfrak g$
\begin{equation}\label{e00}
\omega\, \in\, H^0(E_H,\, \Omega^1_{E_H}\otimes_{\mathbb C}
{\mathfrak g})
\end{equation}
satisfying the following three conditions:
\begin{itemize}
\item $\omega(z)(\zeta_v(z)) \, =\, v$ for all
$v\, \in\, {\mathfrak h}$ and $z\, \in\, E_H$, where $\zeta_v$
is the vector field defined above,

\item $\tau^*_h\omega\, =\, \text{Ad}(h^{-1})\circ\omega$
for all $h\, \in\,H$, where $\tau_h$ is defined in
\eqref{e2}, and

\item for each point $z\, \in\, E_H$, the homomorphism from
the holomorphic tangent space
\begin{equation}\label{e-1}
\omega(z) : T_zE_H\, \longrightarrow\, {\mathfrak g}
\end{equation}
is an isomorphism of vector spaces.
\end{itemize}

All Chern classes will have coefficients in $\mathbb Q$.

By a Calabi--Yau manifold we mean a compact 
K\"ahler manifold $M$ with $c_1(M)\, =\, 0$.

Sorin Dumitrescu \cite{Du} proved that if a Calabi--Yau manifold
has a rigid holomorphic geometric structure of affine algebraic 
type, then it admits a finite covering by a torus.
In \cite{Mc} the second author conjectured a more general result 
that if a Calabi--Yau manifold $M$ admits a holomorphic Cartan
geometry, then $M$ is covered by a torus. Our aim here is to
prove this conjecture under the assumption that $M$ is
complex projective.

This material is based upon works supported by the Science Foundation
Ireland under Grant No. MATF634.

\section{Flat vector bundles on a Calabi--Yau manifold}

Let $M$ be a Calabi--Yau manifold. We will later impose the
condition that $M$ is projective.

\begin{lemma}\label{lem1}
Let $(E\, ,\nabla)$ be a flat vector bundle over 
$M$. Let
$$
E\, \longrightarrow\, Q \, \longrightarrow\, 0
$$
be a quotient holomorphic vector bundle such that $c_1(Q)\,=\, 0$.
Then
$$
c_i(Q)\,=\, 0
$$
for all $i\, \geq\,1$.
\end{lemma}

\begin{proof}
{}It is well known \cite[p. 764, Th\'eor\`eme 2(2)]{Be} that
there is a finite \'etale Galois covering
\begin{equation}\label{e1}
\widetilde{M}\,\longrightarrow\, M
\end{equation}
by a product $\widetilde{M}\,=\, M_{CY} \times T^{n}$ of a 
simply connected Calabi--Yau manifold with a complex
torus. (Note that $X_j$ in \cite[Th\'eor\'eme 2(2)]{Be} are simply
connected \cite[p. 763, Proposition 4(2)]{Be}.)
Clearly we can replace $M$ by $\widetilde{M}$ without
loss of generality, and so assume that the fundamental group 
$\pi_1(M)$ is $\mathbb{Z}^n$. The flat vector bundle
$E$ is given by a linear representation of $\mathbb{Z}^n$.
Since $\mathbb{Z}^n$ is abelian, 
any complex linear representation of $\mathbb{Z}^n$
admits a filtration such that each successive quotient is
a complex one--dimensional $\mathbb{Z}^n$--module. Therefore,
we have a filtration of $E$ by holomorphic subbundles
\begin{equation}\label{e3}
0\,=\, E_0\, \subset\, E_1 \, \subset\, \cdots \, \subset\,
E_{\ell-1} \, \subset\, E_\ell \,=\, E
\end{equation}
such that each successive quotient $E_i/E_{i-1}$, $i\,\in\,
[1\, , \ell]$, is a flat line bundle.

Therefore, $E$ is numerically flat \cite[p. 311,
Theorem 1.18]{DPS}. In particular, $E$ is
numerically effective. Hence
$Q$ is numerically effective \cite[P. 308, Proposition
1.15(i)]{DPS}. Since $c_1(Q)\,=\, 0$, this implies
that $E$ is numerically flat (see 
\cite[p. 311, Definition 1.17]{DPS}). Consequently,
$c_i(Q)\,=\, 0$
for all $i\, \geq\,1$ \cite[p. 311, Corollary 1.19]{DPS}.
\end{proof}

Fix any K\"ahler class on $M$.
By Yau's proof of Calabi's conjecture \cite{Ya}, $M$ admits
a Ricci flat K\"ahler metric, say $\theta$, in the
given K\"ahler class. Note that the
condition that $\theta$ is Ricci flat implies that
the holomorphic tangent bundle $TM$ is polystable of
degree zero with respect to $\theta$ (i.e., a direct
sum of stable vector bundles of degree
zero).

\begin{lemma}\label{lem2}
Let $M$ be a Calabi--Yau complex projective manifold.
Let $E\, \longrightarrow\, M$ be a holomorphic vector bundle with
a holomorphic connection.
Then $E$ admits a holomorphic flat connection.
\end{lemma}

\begin{proof}
This was proved in \cite{Bi} (see \cite[p. 2827,
Theorem A(1)]{Bi}). Note that $TM$ is semistable because it
is polystable. Therefore, \cite[Theorem A(1)]{Bi} applies.
\end{proof}

\section{Cartan geometries and Calabi--Yau manifolds}

As before, $M$ is a Calabi--Yau projective manifold. Let $G$ be a 
complex
Lie group, and let $H\, \subset\, G$ be a closed complex subgroup.
Let $(E_H\, ,\omega)$ be a pair as in \eqref{e0} and \eqref{e00}
defining a Cartan geometry of type $G/H$ on $M$.

Let
\begin{equation}\label{e5}
E_G\, :=\, E_H\times_H G\, \longrightarrow\, M
\end{equation}
be the holomorphic principal $G$--bundle obtained by extending the
structure group of the principal $H$--bundle $E_H$ using the
inclusion map of $H$ in $G$. The form $\omega$ defines a holomorphic
connection on $E_G$. We recall the construction of this
holomorphic connection.

Let $\omega_{\text{MC}}\, :\,
TG\, \longrightarrow\, G\times\mathfrak g$
be the $\mathfrak g$--valued Maurer--Cartan
one--form on $G$ constructed using
the left invariant vector fields.
Consider the $\mathfrak g$--valued holomorphic
one--form
$$
\widetilde{\omega}\, :=\, p^*_1 \omega + p^*_2\omega_{\text{MC}}
$$
on $E_H\times G$, where $p_1$ (respectively,
$p_2$) is the projection of $E_H\times G$ to
$E_H$ (respectively, $G$). This form $\widetilde{\omega}$
descends to a $\mathfrak g$--valued holomorphic 1--form on the
quotient space $E_G$ in \eqref{e5}. This descended 1--form
defines a holomorphic connection on $E_G$.

\begin{theorem}\label{thm1}
Any Calabi--Yau complex projective manifold $M$ bearing a
holomorphic Cartan geometry is holomorphically covered by an 
abelian variety.
\end{theorem}

\begin{proof}
Let $\text{ad}(E_H)$ and $\text{ad}(E_G)$ be the adjoint vector
bundles of $E_H$ and $E_G$ respectively.
Since the homomorphism in \eqref{e-1} is an isomorphism, it
follows that
\begin{equation}\label{e6}
TM\, =\, \text{ad}(E_G)/\text{ad}(E_H)\, .
\end{equation}

We noted above that $E_G$ is equipped with a holomorphic connection.
A holomorphic connection on $E_G$ defines a holomorphic connection
on $\text{ad}(E_G)$. Hence $\text{ad}(E_G)$ admits a holomorphic 
connection. So from Lemma \ref{lem2} it follows that
$\text{ad}(E_G)$ admits a flat holomorphic connection. In view of
Lemma \ref{lem1}, from \eqref{e6} we conclude that $c_2(TM)\,=\,
0$. Now the proof is completed by Lemma 1 of \cite{Mc}.
\end{proof}


\end{document}